\documentclass[12pt,a4paper,twoside,reqno]{amsart}

\usepackage{amsmath, amssymb, amsthm, enumerate}

\newcommand{\IR}{\mathbb{R}}
\newcommand{\IH}{\mathbb{H}}

\DeclareMathOperator{\Ric}{Ric}
\DeclareMathOperator{\tr}{tr}

\DeclareMathOperator{\dist}{dist}
\DeclareMathOperator{\Isom}{Isom}
\DeclareMathOperator{\diam}{diam}
\DeclareMathOperator{\inj}{inj}
\DeclareMathOperator{\Sym}{Sym}
\DeclareMathOperator{\vol}{vol}
\DeclareMathOperator{\const}{const}

\DeclareMathOperator{\DIV}{div}
\newcommand{\dd}{\textup{d}}
\DeclareMathOperator{\diag}{diag}
\newcommand{\EMPTY}[1]{}
\newtheorem{Theorem}{Theorem}[section]
\newtheorem{Lemma}[Theorem]{Lemma}
\newtheorem{Corollary}[Theorem]{Corollary}
\newtheorem{Proposition}[Theorem]{Proposition}
\numberwithin{equation}{section}

\title[Generalized Dehn filling]{Construction of Einstein metrics \\ by generalized Dehn filling}
\author{Richard H Bamler}
\address{Department of Mathematics, Princeton University, Princeton, New Jersey 08544}
\email{rbamler@math.princeton.edu}
\date{November 10, 2009}

\begin{document}

\begin{abstract}
In this paper, we present a new approach to the construction of Einstein metrics by a generalization of Thurston's Dehn filling.
In particular in dimension $3$, we will obtain an analytic proof of Thurston's result.
\end{abstract}
\maketitle

\section{Introduction}
The goal of this paper is to give an analytic construction of Riemannian metrics $g$ which satisfy the Einstein equation $\Ric_g = -(n-1)g$, by a process similar to Thurston's Dehn filling (see \cite{Thu}).
We will first describe the topology of the manifolds on which these metrics live:

Let $(M^n_{hyp},g_{hyp})$ be a hyperbolic manifold of dimension $n \geq 3$ and finite volume.
Denote its cusps by $N_1, \ldots, N_p$ and assume that these are diffeomorphic to $[0,\infty) \times T^{n-1}$, i.e. that the cusps are standard.
We can always choose the $N_k$ so that they are bounded by tori $T_k = \partial N_k$ which are images of horospheres under the universal covering projection and on which $\inj = \mu_n$ where $\inj$ is the injectivity radius and $\mu_n$ the Margulis constant. (For a more detailed description see subsection \ref{sec:Hypmfs}.) 
Now apply the following surgery procedure:
Cut $M_{hyp}$ along the $T_k$, throw away the cusps $N_k$ and glue in $p$ solid tori $\approx D^2 \times T^{n-2}$ by identifying their boundary with the $T_k$.
The topology of the resulting manifold can be uniquely characterized by the homotopy classes of meridional loops $\sigma_k \subset \partial (D^2 \times T^{n-2})$ inside the $T_k$ (i.e. images of loops $S^1 \times \{ pt \} \subset D^2 \times T^{n-2}$ under the gluing identification $\partial (D^2 \times T^{n-2}) \to T_k$).
These homotopy classes are simple, i.e. not a nontrivial multiple of another homotopy class.
Vice versa, given a homotopy class of a simple closed loop $\sigma_k \subset T_k$ for each $k$, we can produce a manifold $M_{\overline{\sigma}} = M_{(\sigma_1, \ldots, \sigma_{p})}$ by this gluing.
In the following we will always assume that the $\sigma_1, \ldots, \sigma_{p}$ are geodesic representatives (inside $T_1, \ldots, T_{p}$) of their homotopy classes and set $\ell_k := \ell(\sigma_k)$ and $\ell_{\min \text{resp.} \max} := (\min \text{resp.} \max)(\ell_k)$.

The statement of the theorem which we are going to prove, is now:
\begin{Theorem} \label{Thm:main}
There is a constant $L = L(n, V)$ such that whenever $\vol M_{hyp} < V$ and $\ell_{\min} > L$, the manifold $M_{\overline{\sigma}}$ carries an Einstein metric $g_{\overline{\sigma}}$.

Moreover, the metrics on $M_{\overline{\sigma}}$ can be constructed in such a way that as $\ell_{\min} \to \infty$, the $(M_{\overline{\sigma}}, g_{\overline{\sigma}})$ converge to the initial hyperbolic manifold $(M_{hyp}, g_{hyp})$ in the pointed Gromov-Hausdorff sense if the basepoints are chosen away from the cusps.
\end{Theorem}

A slightly weaker statement was also claimed in \cite{And}.
Theorem \ref{Thm:main} immediately implies the Dehn filling Theorem in dimension $3$:
\begin{Corollary}
Let the dimension $n=3$.
There is a constant $L=L(V)$ such that whenever $\vol M_{hyp} < V$ and $\ell_{\min} > L$, the manifold $M_{\overline{\sigma}}$ is hyperbolic.
\end{Corollary}

A slightly weaker version of this theorem was proven by Thurston (\cite{Thu}) using the deformation theory of Kleinian groups.
Our methods provide a new and analytic proof of his result.

We will give a short sketch of the proof of Theorem \ref{Thm:main}:
First, we endow the solid tori which we will glue into the hyperbolic manifold $M_{hyp}$ with a special Einstein metric called the \emph{black-hole metric.}
This metric is asymptotically hyperbolic to its end and thus each gluing can be arranged to be arbitrarily smooth for large $\ell_{k}$.
Hence, the resulting metric is almost Einstein, i.e. its traceless Ricci tensor is small in some $C^{m, \alpha}$-sense.
Eventually, we apply an inverse function theorem like argument to perturb the metric into the desired Einstein metric.

We mention that our proof builds on previous work of Tian (\cite{Tia}) and Anderson (\cite{And}).
Tian established the $3$-dimensional case in which $M_{hyp}$ has only one cusp.
Later Anderson described a construction for the higher dimensional case and developed new analytical tools of which we will also partly make use here.

We want to point out that the case in which the hyperbolic manifold $M_{hyp}$ has more than one cusp, is substantially more difficult than the case of one cusp for the following reason:
The accuracy of the gluing in the first step (i.e. the construction of the almost Einstein metric) depends polynomially on the minimum $\ell_{\min}$ of the $\ell_k$.
However, as the $\ell_k$ get large, the invertibility of the linearized Einstein equation deteriorates logarithmically in the maximum $\ell_{\max}$ of the $\ell_k$.
So in the case of one cusp $\ell_{\min} = \ell_{\max}$ and thus the accuracy of the gluing increases more rapidly than the invertibility deteriorates.
But if $M_{hyp}$ has more than one cusp and $\ell_{\min}$, $\ell_{\max}$ are not sufficiently controlled with respect to one another, then this consideration fails.
In \cite{And}, Anderson sketches an argument how to get around this issue by looking at certain moduli spaces of solutions of a modified Einstein equation.
In this paper, we will be able to deal with the problems that arise in this multiple cusp case and we will give a complete proof of Theorem \ref{Thm:main}.
In fact, our argument will be more elementary and we find it a more natural way of looking at the problem.

The idea behind our proof is that the reason for the bad invertibility of the linearized Einstein equation lies in certain variations of the metric (so called \emph{trivial Einstein variations}) which correspond to a change of the moduli of the cross-sectional tori of the cusps.
It will turn out that with respect to some cleverly chosen norms (see section \ref{sec:unifnorms}), which treat these trivial Einstein variations separately, the invertibility of the linearized Einstein equation becomes in fact independent of $\ell_{\max}$ (see our Proposition \ref{Prop:Linv} as opposed to Proposition 3.2 in \cite{And}).
However, these new norms make it necessary to reprove the inverse function theorem in order to be applicable to our setting (see section \ref{sec:Applift}).

Another difference between our and Anderson's proof is that we have replaced the proof of Lemma 3.4 in \cite{And} which used to involve the theory of the moduli spaces of conformally compact Einstein metrics, by an elementary argument by which we can even show a slight generalization.

We remark that it still remains an interesting question whether the constant $L$ in Theorem \ref{Thm:main} can be chosen independent of the volume of $M_{hyp}$.
Hodgson and Kerckhoff (see \cite{HK}) could confirm this in dimension $3$ using algebraic techniques.

The paper is organized as follows:
In section \ref{sec:Prel}, we will review some basic facts which were also used in \cite{And}.
Section \ref{sec:constrprocess} contains a brief recapitulation of the construction of the almost Einstein metric as described in \cite{And}.
In sections \ref{sec:unifnorms} to \ref{sec:estforLm} we carry out the main argument.
In order to keep these chapters concise, we will defer most of the technical calculations to sections \ref{sec:notoncusp} and \ref{sec:notonBH}.

I want to thank my advisor Gang Tian for calling my attention on this subject and for his continual support.
I am also indebted to the HIM in Bonn for their hospitality and Hans-Joachim Hein for many inspiring discussions.

\section{Preliminaries} \label{sec:Prel}
\subsection{Hyperbolic manifolds} \label{sec:Hypmfs}
We recall the \emph{thick-thin-decomposition} for hyperbolic manifolds
\begin{Theorem}
There is a constant $\mu_n > 0$, the \emph{Margulis constant,} such that the following holds:
If $M^n_{hyp}$ is a finite volume hyperbolic manifold then $M_{hyp}$ can be decomposed into a \emph{thin part} $M_{thin}$ and a \emph{thick part} $M_{thick}$ with $M_{hyp} = M_{thin} \dot \cup M_{thick}$ such that:
\begin{itemize}
 \item $\inj \geq \mu_n$ on $M_{thick}$ and $M_{thick}$ is relatively compact in $M_{hyp}$.
 \item $M_{thin}$ is a finite union of connected open sets $N_1, \ldots, N_p$ and $N'_1, \ldots, N'_{p'}$ where
 \begin{itemize}
 \item  the $N_k$ are cusps of the form $[0, \infty) \times (T^{n-1} / \Gamma_k)$ for finite subgroups $\Gamma_k < \Isom T^{n-1}$ with a warped product metric
 \begin{equation} \label{eq:warpedproductrep} g_{hyp} = ds^2 + e^{-2s} g_{flat, T^{n-1} / \Gamma_k}. \end{equation}
 In the case in which $\Gamma_k = \{ 1 \}$, we call $N_k$ \emph{standard}.
\item and the $N'_k$ are covered by cylindrical neighborhoods around geodesics in hyperbolic space.
\end{itemize}
Furthermore, we can choose the $N_k$ such that their boundaries are images of horospheres under the universal covering projection and such that $\inj = \mu_n$ at some point on $\partial N_k$.
In the case in which $N_k$ is standard this implies that $\inj = \mu_n$ on $\partial N_k$.
\end{itemize}
In every dimension, $\diam M_{thick}$ is bounded from above by a constant which only depends on an upper bound on $\vol M_{hyp}$ and in dimension $n \geq 4$, this is even true for the diameter of $M_{thick} \cup N'_1 \cup \ldots \cup N'_{p'} = M_{hyp} \setminus N_1 \cup \ldots \cup N_p$.
\end{Theorem}

We can compute the volume of the cusps in terms of the area of their boundary surface:
\begin{Lemma}
There is a constant $\eta_n$ such that for all cusps $N_k$ we have
\[ \vol N_k = \eta_n \vol \partial N_k. \]
\end{Lemma}
\begin{proof}
This can be checked easily using (\ref{eq:warpedproductrep}).
\end{proof}

So a bound on the volume of $M_{hyp}$ gives us a bound on the volume of the $\partial N_k$.
Since $\partial N_k$ lies in the thick part, we have a bound on the injectivity radius of $\partial N_k$ (which is slightly larger than $\mu_n$ since $\partial N_k$ is not totally geodesic). 
The next lemma shows that in fact we get a bound on the diameter of $\partial N_k$ from an upper volume bound on $M_{hyp}$.
This implies furthermore, that the tori that can occur as cusp cross-sections of a hyperbolic manifold with a given volume bound form a bounded subset in the moduli space of flat tori.
\begin{Lemma}
For every $V < \infty$ and $\iota > 0$ there is a $d (n, V, \iota) < \infty$ such that for any flat torus $T^{n-1}$ we have
\[ \vol T^{n-1} < V \quad \text{and} \quad \inj T^{n-1} > \iota \quad \Longrightarrow \quad \diam T^{n-1} < d. \]
\end{Lemma}
\begin{proof}
Let $\gamma : [0,l] \to T^{n-1}$ be a minimizing geodesic.
Then the balls $B_{\iota}(\gamma(\iota))$, $B_{\iota}(\gamma(3 \iota))$, \ldots are pairwise disjoint and have volume $\omega_{n-1} \iota^{n-1}$.
So $l < 2(\frac{V}{\omega_{n-1} \iota^{n-1}} + 1) \iota$.
\end{proof}

\subsection{The Einstein operator} \label{subsec:Einstop}
For any symmetric bilinear form $h$ and any 1-form $\alpha$ on a Riemannian manifold $(M,g)$ we define the divergence and its formal conjugate by ($e_k$ is a local orthonormal frame field)
\[ \DIV_g (h) = - \sum_{i=1}^n (\nabla_{e_i} h)(e_i, \cdot), \qquad (\DIV_g^* \alpha)(X,Y) = \tfrac12 ( (\nabla_X \alpha)(Y) + (\nabla_Y \alpha)(X) ). \]
Observe that $\DIV_g^* \alpha = \frac12 \mathcal{L}_{\alpha^\sharp} g$.
Let $h$ be a bilinear form.
We can express the derivative of the Ricci curvature in the direction of $h$ by (for a computation see \cite[sec 2.3]{Top})
\[ \dd\Ric_g(h) = - \tfrac12 \triangle_L h - \DIV_g^* ( \DIV_g h + \tfrac12 \dd \tr_g h ). \]
Here $(\triangle_L h)(X,Y) = (\triangle h)(X,Y) + 2 R(h)(X,Y) - h(\Ric(X),Y) - h(X, \Ric(Y))$ is the Lichnerowicz Laplacian and $R(h)(X,Y) = \tr h(R(\cdot, X)Y, \cdot)$.
Since computing the Ricci tensor is a diffeomorphism invariant operation, we have for any 1-form $\alpha$
\begin{equation} \dd\Ric_g(\DIV_g^* \alpha) = \tfrac12 \mathcal{L}_{\alpha^{\sharp}} \Ric_g. \label{eq:Ricdiffeo} \end{equation}
Thus $\dd \Ric_g$ is not an elliptic operator.
In order to make it elliptic, we have to add an extra term:
Let $\overline{g}$ be an arbitrary fixed background metric on $M$.
We define $\Psi_{\overline{g}} : \{ g \in C^\infty (M ; \Sym_2 T^*) \; : \; g>0 \} \to C^\infty (M ; \Sym_2 T^*)$ by
\[ \Psi_{\overline{g}}(g) := \Ric_g + \DIV_g^*(\DIV_{\overline{g}} g + \tfrac12 \dd \tr_{\overline{g}} g). \]
Its derivative at $\overline{g}$ is
\[ (\dd\Psi_{\overline{g}})_{\overline{g}}(h) = - \tfrac12 \triangle_L h \]
hence elliptic.
For our purposes we define the Einstein operator $\Phi_{\overline{g}} : \{ g \in C^\infty (M ; \Sym_2 T^*) \; : \; g>0 \} \to C^\infty (M ; \Sym_2 T^*)$ by
\[ \Phi_{\overline{g}}(g) = \Psi_{\overline{g}}(g) + (n-1) g. \]
We have
\begin{multline*} 
(\dd\Phi_{\overline{g}})_{\overline{g}}(h) = - \tfrac12 \triangle_L h + (n-1)h \\
= \tfrac12 \left(- \triangle h - 2 R(h) + \Ric \circ h + h \circ \Ric + 2(n-1)h \right).
\end{multline*}
Set $L_{\overline{g}} := 2(\dd\Phi_{\overline{g}})_{\overline{g}}$ and call elements in the kernel of $L_{\overline{g}}$ \emph{Einstein variations.}
Using a Weitzenb\"ock formula, we can express this linear operator as
\[ L_{\overline{g}} h = (\DIV^* \DIV + d^* d)h - R(h) + \tfrac12 \Ric \circ h + \tfrac12 h \circ \Ric + 2(n-1) h \]
where $d : C^\infty(M; \Sym_2 T^*) \to C^\infty(M; \Lambda_2 T^* \otimes T^*)$ and its formal conjugate $d^* : C^\infty(M; \Lambda_2 T^* \otimes T^*) \to C^\infty(M; \Sym_2 T^*)$ are defined by
\begin{multline*} 
 (dh)(X,Y,Z) = (\nabla_X h)(Y, Z) - (\nabla_Y h)(X, Z),\\
 (d^*t)(X, Y) = - \tfrac12 \sum_{i=1}^n ( (\nabla_{e_i} t)(e_i, X, Y) + (\nabla_{e_i} t)(e_i, Y, X) ). 
\end{multline*}
If $\overline{g}$ is Einstein with $\Ric_{\overline{g}} = -(n-1) \overline{g}$, we have
\[ L_{\overline{g}} h = - \triangle h - 2 R(h) = (\DIV^* \DIV + d^* d)h - R(h) + (n-1) h. \]
Tracing this equation gives us
\begin{equation} \tr L_{\overline{g}} h = \nabla^*\nabla \tr h + 2(n-1) \tr h. \label{eq:tr} \end{equation}
If $\overline{g}$ is hyperbolic of constant sectional curvature $-1$, we get
\begin{equation} L_{\overline{g}} h = - \triangle h - 2 h + 2 (\tr_{\overline{g}} h) \overline{g} 
= (\DIV^* \DIV + d^* d)h + (\tr h)_{\overline{g}} \overline{g} + (n-2) h. \label{eq:hypWbck} \end{equation}

\begin{Lemma} \label{Lem:dphiEinst}
If $M$ is closed and $\Ric_g < 0$, then $\Phi_{\overline{g}}(g) = 0$ implies $\Ric_g = -(n-1)g$.
\end{Lemma}
\begin{proof}
This Lemma can be found in \cite{Biqbook} and \cite[Lemma 2.1]{And}.
 We copy the proof from the latter source since we need a variation of the argument later on.
 Let $\beta_{g} (h) := \DIV_{g} h + \frac12 \dd \tr_{g} h$ be the Bianchi operator.
 Applying $\beta_g$ to $\Phi_{\overline{g}}(g) = 0$ yields
 \begin{equation} 
 0 = \beta_{g}(\Phi_{\overline{g}}(g)) = \beta_g \DIV_g^* \beta_{\overline{g}} (g) = \tfrac12 \left( \nabla^*\nabla \beta_{\overline{g}}(g) - \Ric_g (\beta_{\overline{g}}(g)) \right). \label{eq:beta}
 \end{equation}
 So $\beta_{\overline{g}}(g) = 0$ and the claim follows.
\end{proof}
Thus, in order to construct Einstein metrics, it suffices to look for zeros of $\Phi_{\overline{g}}$.

A similar result to Lemma \ref{Lem:dphiEinst} holds in the differential sense:
\begin{Lemma} \label{Lem:Einstvar}
If $(M, \overline{g})$ is a complete Einstein manifold and $h$ a symmetric bilinear form such that $|h|(x), |\nabla h|(x) \to 0$ for $x \to \infty$, then $L_{\overline{g}} h = 0$ implies
\[ \dd \Ric_{\overline{g}} (h) = -(n-1) h, \qquad \DIV_{\overline{g}} h = 0, \qquad \tr_{\overline{g}} h = 0. \]
\end{Lemma}
\begin{proof}
 The proof is the same as in \cite[Lemma 3.6]{AndCC}.
 Differentiate (\ref{eq:beta}) with respect to $g$ to find
 \begin{equation} 0 = \beta_{\overline{g}} L_{\overline{g}} h = \nabla^* \nabla \beta_{\overline{g}} h + (n-1) \beta_{\overline{g}} h. \label{eq:bg} \end{equation}
 So $\beta_{\overline{g}} h = 0$.
 Moreover, by (\ref{eq:tr}) we conclude $\tr h = 0$.
\end{proof}

Observe that conversely not every Einstein variation is divergence or trace free.

\subsection{The hyperbolic cusp}\label{subsec:hypcusp}
We introduce a representation for the metric $g_{hyp}$ on a hyperbolic cusp which is different from (\ref{eq:warpedproductrep}):
Consider the coordinates $(r,x_2, \ldots, x_n)$ on $\IR_+ \times \IR^{n-1}$ and the hyperbolic metric
\begin{equation} \label{eq:hypmetricrep} g_{hyp} = r^{-2} \dd r^2 + r^2 (\dd x_2^2 + \ldots + \dd x_n^2).
\end{equation}
Note that in these coordinates, $g_{hyp}$ is not conformally equivalent to the Euclidean metric (as opposed to the coordinates that arise after the transformation $r \to \frac1r$).
Obviously, the metric is invariant under the action by Euclidean isometries on the last factor.
Now every hyperbolic cusp is the quotient of $\IR_+ \times \IR^{n-1}$ metric under a discrete subgroup of those isometries.

We will be interested in Einstein deformations of the metric $g_{hyp}$ which are invariant by the group $\IR^{n-1}$ of Euclidean tranformations on the last factor.
One type of deformation will be very essential:
Let $u_{ij}$ be a symmetric $(n-1) \times (n-1)$ matrix indexed by $i, j = 2, \ldots, n$.
Then if $u_{ij} > -\delta_{ij}$
\[ r^{-2} \dd r^2 + r^2 (\dd x_2^2 + \ldots + \dd x_n^2 + u_{ij} \dd x_i \dd x_j) \]
is isometric to $g_{hyp}$ hence it is also Einstein.
We will denote this metric by $g_{hyp} + u$.
It can be checked (e.g. using (\ref{eq:tr})) that the equation $\Phi_{g_{hyp}}(g_{hyp} + u) = 0$ is equivalent to $\tr u = 0$.
Likewise, dropping the lower bound for $u_{ij}$ and setting $h = r^2 u_{ij} \dd x_i \dd x_j$, we find that $L_{g_{hyp}} h = 0$ iff $\tr u = 0$.
We will call variations of this kind \emph{trivial Einstein variations.}

\subsection{The black-hole metric} \label{subsec:BH}
We recall the definition of the black-hole metric $(M_{BH}, g_{BH})$ as in \cite{And}.
Introduce coordinates $(r, \theta, x_3, \ldots, x_n)$ on $\IR^n = \IR^2 \times \IR^{n-2}$.
Here, $(r, \theta)$ denote polar coordinates on the first factor such that $r$ is running from $r_+ = 2^{1/n-1}$ to $\infty$ and $\theta$ from $0$ to $\beta = \frac{4 \pi}{(n-1) r_+}$.
The black-hole metric is defined as
\[ g_{BH} = V^{-1} \dd r^2 + V \dd \theta^2 + r^2 (\dd x_3^2 + \ldots + \dd x_n^2)
\]
where
\[ V(r) = r^2 - \frac{2}{r^{n-3}}. \]
Using the coordinate transformation $r = r_+ + \frac12 s^2$, we obtain
\[ g_{BH} = \frac{s^2}{V(r_+ + \frac12 s^2)} \dd s^2 + s^2 \frac{V(r_+ + \frac12 s^2)}{s^2} \dd \theta^2 + (r_+ + \tfrac12 s^2)^2 (\dd x_3^2 + \ldots + \dd x_n^2). \]
Since $\frac{V(r_+ + \frac12 s^2)}{s^2} \to \frac12 (n-1) r_+$ as $s \to 0$, we conclude that $g_{BH}$ is in fact smooth at the origin.
The sectional curvatures of $g_{BH}$ are
\begin{alignat*}{2}
 K_{12} &= -1 + \frac{(n-3)(n-2)}{r^{n-1}}, \qquad & K_{1i} &= K_{2i} = -1 - \frac{n-3}{r^{n-1}}, \qquad i \geq 3, \\
 K_{ij} &= -1 + \frac{2}{r^{n-1}}, \qquad i,j \geq 3 &&
\end{alignat*}
Furthermore, $g_{BH}$ is Einstein with $\Ric_{g_{BH}} = -(n-1) g_{BH}$.
For $n=3$ this metric is just the standard hyperbolic metric in cylindrical coordinates.

Observe that away from the origin the metric $g_{BH}$ is asymptotic to the standard hyperbolic metric $g_{hyp}$ from subsection \ref{subsec:hypcusp}.
To be precise: The black-hole manifold minus a large cylinder around the core $\IR^{n-2}$ is geometrically close to some subset of $\IH^n / \langle \gamma \rangle$ where $\gamma$ is a parabolic transformation.
Taking the hyperbolic metric $g_{hyp}$ as a background metric and identifying $\theta$ with $x_2$, we find
\begin{multline*}
 g_{BH} - g_{hyp} = (V^{-1}(r) - r^{-2}) \dd r^2 + (V(r) - r^2 ) \dd x_2^2 \\ = 2 r^{-n+1} \frac{r^2}{V(r)} \cdot r^{-2} \dd r^2 - 2 r^{-n+1} \cdot r^2 \dd x_2^2.
\end{multline*}
From this, we can conclude $| \nabla^m(g_{BH} - g_{hyp} ) | = O(r^{-n+1})$ for all $m \geq 0$.

We analyze the behaviour of this metric under the addition of small trivial Einstein deformations.
Let $u_{ij}$ ($i,j = 2, \ldots, n$) be a traceless symmetric $(n-1) \times (n-1)$ matrix and set
\[ g_{BH} + u = g_{BH} + r^2 u_{ij} \dd x_i \dd x_j \]
where we set $x_2 = \theta$.
Note that $g_{BH} + u$ is only smooth away from the origin.
By the closeness of $g_{BH}$ to $g_{hyp}$ we find for sufficiently small $u$ and say $r > r_+ + 1$
\begin{equation} \Phi_{g_{BH}} (g_{BH} + u) = |u| O(r^{-n+1}) \qquad \text{and} \qquad L_{g_{BH}} u = |u| O(r^{-n+1}). \label{eq:trivEinsdecay} \end{equation}
Since we will need it later, we mention the following bound:
Let $u$ be small and $u'$ be another traceless symmetric $(n-1) \times (n-1)$ matrix.
Then
\begin{equation} \left| \left( \dd \Phi_{g_{BH}} \right)_{g_{BH} + u} (u') - \left( \dd \Phi_{g_{BH}} \right)_{g_{BH}} (u') \right| = |u| |u'| O(r^{-n+1}). \label{eq:trivEinsnoinfl}
\end{equation}
The same decay holds for all higher covariant derivatives of the left hand side.

It will be useful later to discuss the geometric quotients of the black-hole metric.
In dimension $n=3$ the isometries of $(M_{BH}, g_{BH})$ are the isometries of hyperbolic space.
For $n > 3$, all isometries leave the topological splitting $M_{BH} = \IR^2 \times \IR^{n-2}$ and the coordinate $r$ invariant, so they act as a rotation or reflection in the origin on the $\IR^2$ factor and as a Euclidean translation on the $\IR^{n-2}$-factor.

Let $T^{n-1}$ be an arbitrary flat torus and $\sigma \subset T^{n-1}$ a simple closed geodesic.
Given this data, we will now construct a Riemannian manifold with boundary $\widehat{N} = \widehat N_{T^{n-1}, \sigma}$ such that
\begin{enumerate}[(i)]
\item $\widehat{N}$ is diffeomorphic to a solid torus $D^2 \times T^{n-2}$.
\item The boundary $\partial \widehat{N}$ is isometric to the given $T^{n-1}$.
\item Denote by $M_{BH}(r \leq R)$ the part of $M_{BH}$ on which $r \leq R$.
Then there is a group of isometries $\Gamma$ of $M_{BH}$ and a number $R > 0$ such that $\widehat{N} = M_{BH}(r\leq R) / \Gamma$.
\item All isometries of $\Gamma$ keep the splitting $M_{BH}  = \IR^2 \times \IR^{n-2}$ and the coordinate $r$ invariant, hence $T^{n-1} \cong \partial \widehat N = M_{BH}  (r = R) / \Gamma$.
\item There is an isometric identification of $M_{BH} (r = R)/ \Gamma$ with $T^{n-1}$ that sends the meridional loop $M_{BH} (r= R, x_3 = \ldots = x_n = 0) / \Gamma$ to $\sigma$.
\end{enumerate}
We will later use $\widehat{N}$ to fill in the truncated cusps.

The manifold $\widehat{N}$ is constructed as follows:
Since $V : [r_+, \infty) \to [0, \infty)$ is invertible, we can find an $R > 0$ such that $V(R) = \left( \ell(\sigma) / \beta \right)^2$. 
Then $M_{BH}(r = R) \cong S^1(\ell(\sigma)) \times \IR^{n-2}$.
Consider the cyclic subgroup $\langle [\sigma] \rangle < \pi_1 (T^{n-1})$ and denote by $\widetilde{T} \to T^{n-1}$ the corresponding cover.
Obviously, $\sigma \subset T^{n-1}$ can be lifted to a closed geodesic loop $\widetilde{\sigma} \subset \widetilde{T}$.
From this we conclude that there is an isometric identification of $M_{BH}(r = R)$ with $\widetilde{T}$ such that the loop $M_{BH}(r=R, x_3 = \ldots = x_n)$ is sent to $\widetilde{\sigma}$.
Consider now the group of deck transformations $\Gamma$ of $\widetilde{T} \to T^{n-1}$.
Its action on $M_{BH}(r = R) \cong \widetilde{T}$ can be uniquely continued to an isometric action on $M_{BH}$ and by the lifting property of $\sigma$ we know that this continuation is even fixed point free.
Hence $M_{BH}/\Gamma$ is smooth, and the manifold $\widehat N := M_{BH} (r \leq R)/ \Gamma$ satisfies the conditions (i)-(v) above.

We remark that the image of the core $\IR^{n-2}$ under the quotient map $M_{BH} \to M_{BH} / \Gamma$ is a torus $\widehat T^{n-2} = r^{-1}(r_+)$ which we call the \emph{core torus.}
Furthermore, all level sets $r^{-1}(r')$ for $r' > r_+$ are diffeomorphic to $T^{n-1}$ and we can check that $(r')^{-1} \diam r^{-1}(r')$ is an increasing function in $r'$ (here $\diam$ denotes the intrinsic diameter).

\section{The construction process} \label{sec:constrprocess}
We will briefly explain how the approximate Einstein metric on the manifold $M_{\overline{\sigma}}$ is constructed.
Recall that we are given simple closed geodesics $\sigma_k$ inside the tori $T_k \subset M_{hyp}$ which bound the cusps $N_k$, and that $\inj = \mu_n$.
In dimension $3$ it is also important to choose the $N'_k$ such that on their boundary tori $T'_k$ we also have $\inj = \mu_n$.

As mentioned in subsection \ref{subsec:BH}, we can find parameters $R_k$ as well as lattices $\Gamma_k < \Isom M_{BH}$ such that the $\widehat N_k := \widehat N_{T_k, \sigma_k} = M_{BH} (r \leq R_k)/ \Gamma_k$ are topological solid tori with boundary isometric to $T_k$ and such that the $\sigma_k$ correspond to meridians.
Set $R_{\min} := \min R_k$ and observe that $R_{\min} \to \infty$ as $\ell_{\min} \to \infty$.

If we glue together the components $M_{hyp} \setminus \bigcup_{k=1}^p N_k$ and $\widehat N_k$, we obtain the manifold $M_{\overline{\sigma}}$.
We can endow $M_{\overline{\sigma}}$ with an almost Einstein metric $g_{\overline{\sigma}}$ in the following sense (to simplify notation, we will denote this metric by $g_{\overline{\sigma}}$ rather than the final Einstein metric):
$g_{\overline{\sigma}}$ equals $g_{hyp}$ on the first component and $g_{BH}$ on the $\widehat N_k$ except on the tubular neighborhoods $\bigcup_{k=1}^p B_1 T_k \cap \widehat N_k$ of radius $1$ around the $T_k$ where an  interpolation between $g_{hyp}$ and $g_{BH}$ is taking place.
Thus $g_{\overline{\sigma}}$ satisfies the Einstein equation on the complement of $\bigcup_{k=1}^p B_1 T_k \cap \widehat N_k$ whereas on $\bigcup_{k=1}^p B_1 T_k \cap \widehat N_k$ the quantity $\Ric_{g_{\overline{\sigma}}} + (n-1) g_{\overline{\sigma}}$ and hence $\Phi_{g_{\overline{\sigma}}}(g_{\overline{\sigma}})$ is very small.
To be precise: for any $m$ there is a $C_m$ such that using the $C^{m,\alpha}$ norms defined in section \ref{sec:unifnorms} below, we have
\[ \Vert \Ric_{g_{\overline{\sigma}}} + (n-1) g_{\overline{\sigma}} \Vert_{m, \alpha} < C_m R_{\min}^{-n+1} \qquad \text{and} \qquad \Vert \Phi_{g_{\overline{\sigma}}} (g_{\overline{\sigma}} ) \Vert_{m, \alpha} < C_m R_{\min}^{-n+1}. \]

Note that in dimension $n=3$, the parts $N'_k$ are already isometric to $\widehat N'_k = \widehat N_{T'_k, \sigma'_k}$ for certain loops $\sigma'_k \subset T'_k$.

For further details of this construction we refer to \cite{And}.

\section{Uniform norms on $M_{\overline{\sigma}}$} \label{sec:unifnorms}
In the following let $L := L_{g_{\overline{\sigma}}}$ and fix some $m \geq 2$ and $0 < \alpha < 1$.
We will further fix an upper volume bound on $M_{hyp}$ and call all constants \emph{uniform} which only depend on this bound, but not on $M_{hyp}$ or $\overline{\sigma}$.
Observe that the Riemannian manifolds $(M_{\overline{\sigma}}, g_{\overline{\sigma}})$, as constructed in the last section, satisfy the following uniform geometric bounds: The conjugate radius is uniformly bounded from below by some positive constant $2\zeta$ and there are uniform bounds $C_m$ such that $\Vert \nabla^m R \Vert < C_m$.

Let $h$ be a symmetric bilinear form on $M_{\overline{\sigma}}$ and $x \in M_{\overline{\sigma}}$.
Pull back the bundle $\textnormal{Sym}_2 T^*$ and its section $h$ to the universal cover $\widetilde{B}_{\zeta}(x)$ of $B_\zeta(x)$.
Choose exponential coordinates on  $\widetilde{B}_{\zeta}(x)$ and trivialize $\textnormal{Sym}_2 T^*$ by parallel transport.
We can now view $h$ as a vector-valued function on a ball $B_\zeta(0) \subset \IR^n$.
Fix the constant $\zeta$ once and for all and define the local H\"older (semi)-norm of $h$ at $x$ by this representation:
\[ \Vert h \Vert_{m, \alpha; x} := \Vert h |_{B_\zeta(x)} \Vert_{m, \alpha}. \]
We note that we have Schauder estimates for these semi-norms:
\begin{equation} \Vert h \Vert_{m, \alpha; x} \leq C \sup_{x' \in B_{\zeta}(x)} \left( \Vert L h \Vert_{m-2, \alpha; x'} + \Vert h \Vert_{0; x'} \right) \label{eq:basicSchauder} \end{equation}
such that $C$ is a uniform constant.

Using these semi-norms it is now easy to define the global H\"older norm by
\[ \Vert h \Vert_{m, \alpha} := \sup_{x \in M_{\overline{\sigma}}} \Vert h \Vert_{m, \alpha; x}. \]

We will need another norm that guarantees a certain decay away from the thick part and the core tori.
We therefore introduce a weight function $W$ (or rather the inverse of a weight function) on $M_{\overline{\sigma}}$ such that for $n>3$
\[ W = \left\{ 
\begin{array}{ll}
( \frac{r}{R_k} )^{0.1} + r^{-0.1} & \text{on $\widehat N_k$} \\
1 & \text{on $M_{\overline{\sigma}} \setminus \bigcup_{k=1}^p \widehat N_k$}
\end{array}
\right.\]
In dimension $3$ we also choose the weight $( \frac{r}{R'_k} )^{0.1} + r^{-0.1}$ on the $N'_k$.
$W$ is not continuous at the $T_k$.
However, this discontinuity will not be essential since the jump is between $1$ and $1+R_k^{-0.1}$ and $R_k > r_+$.
On each $\widehat N_k$ the weight function $W$ attains its minimum at $r = R_k^{1/2}$.
For later use, choose points $c_k \in \widehat N_k( r= R_k^{1/2})$ (and $c'_k \in N'_k(r = (R'_k)^{1/2})$ in dimension $3$).
They lie approximately in the centers of the $\widehat N_k$.
Set
\[ \Vert h \Vert_{m, \alpha; *} := \sup_{x \in M_{\overline{\sigma}}} W^{-1}(x) \Vert h \Vert_{m, \alpha; x}. \]
It is immediate that we can derive uniform Schauder estimates for the norms $\Vert \cdot \Vert_{m, \alpha}$ and $\Vert \cdot \Vert_{m, \alpha; *}$ from (\ref{eq:basicSchauder}):
\begin{alignat*}{1}
 \Vert h \Vert_{m, \alpha} &\leq C \left( \Vert L h \Vert_{m-2, \alpha} + \Vert h \Vert_0 \right) \\
 \Vert h \Vert_{m, \alpha; *} &\leq C \left( \Vert L h \Vert_{m-2, \alpha; *} + \Vert h \Vert_{0; *} \right) 
\end{alignat*}

Finally, we have to define a more complicated norm that guarantees decay towards some trivial Einstein variation:
Let $\rho_1, \ldots, \rho_p$ be cutoff functions on $M_{\overline{\sigma}}$ such that $\rho_k \equiv 1$ on $\widehat N_k \setminus ( B_1 T_k \cup B_2 \widehat T^{n-2}_k )$ and $\rho_k \equiv 0$ on $M_{\overline{\sigma}} \setminus \widehat N_k$ and $B_1 \widehat T^{n-2}_k$ where $\widehat T^{n-2}_k$ is the core torus of $\widehat N_k$.
We may assume that the $\rho_k$ are constructed in such a way that they satisfy some universal $C^m$ bound for each $m$.
Let $u_1, \ldots, u_p$ be trivial Einstein variations of the hyperbolic cusp metric which we assume to be defined on the corresponding $\widehat N_k$.
Represent $h$ by
\begin{equation}
h = \bar h + \sum_{k=1}^p \rho_k u_k. \tag{$*$}
\end{equation}
and define
\[ \Vert h \Vert_{m, \alpha; **} := \inf_{\stackrel{\scriptstyle \bar h, u_1, \ldots, u_p}{\text{satisfy ($*$)}}} \bigg( \Vert \bar h \Vert_{m, \alpha; *} + \sum_{k=1}^p | u_k | \bigg) \]
where we use an arbitrary uniform norm on the (finite dimensional) space of trivial Einstein deformations.
In dimension $3$ we have to alter the definition in order to also consider trivial Einstein variations $u'_k$ on the $N'_k$.
Observe that since $W^{-1} > c > 0$ we have for some uniform $C$
\[ C^{-1} \Vert h \Vert_{m, \alpha} \leq \Vert h \Vert_{m, \alpha; **} \leq  \Vert h \Vert_{m, \alpha; *}. \]

\begin{Lemma}
We have the following uniform Schauder estimate for \hbox{$\Vert \cdot \Vert_{m, \alpha; **}$}:
\begin{equation} \Vert h \Vert_{m, \alpha; **} \leq C \left( \Vert L h \Vert_{m-2, \alpha; *} + \Vert h \Vert_{0; **} \right). \label{eq:ssSchauder} \end{equation}
Note that the second norm is a $*$-norm.
\end{Lemma}
\begin{proof}
 We carry out the proof for $n>3$ (for $n=3$ we have to consider the $N'_k$ as well).
 Choose a decomposition $h = \bar h + \sum_{k=1}^p \rho_k u_k$.
 From (\ref{eq:trivEinsdecay}) we find that $\Vert L \rho_k u_k \Vert_{m-2, \alpha; *} \leq C |u_k|$.
Hence
\[ \Vert L \bar h \Vert_{m-2, \alpha; *} \leq \Vert L h \Vert_{m-2, \alpha;*} + C \sum_{k=1}^p |u_k|. \]
So by the Schauder estimate for $\Vert \cdot \Vert_{m, \alpha; *}$ we find
\[ \Vert \bar h \Vert_{m, \alpha; *} \leq C \bigg( \Vert L h \Vert_{m-2, \alpha; *} + \sum_{k=1}^p |u_k| + \Vert \bar h \Vert_{0;*} \bigg) \]
hence the conclusion.
\end{proof}

The following Lemma gives us a tool to estimate the $\Vert \cdot \Vert_{m, \alpha; **}$ norm:
\begin{Lemma} \label{Lem:trivcenter}
 Let $h$ be a symmetric bilinear form on $M_{\overline{\sigma}}$.
Choose $u_k$ such that $| (h- u_k)(c_k)  |$ is minimal for each $k$ and set $\bar h = h - \sum_{k=1}^p \rho_k u_k$.
In dimension $3$ also consider $u'_k$ such that $| (h-u'_k) (c'_k) |$ is minimal.
Then there is a uniform constant $C$ such that
\[  \Vert h \Vert_{m, \alpha; **} \leq \bigg( \Vert \bar h \Vert_{m, \alpha; *} + \sum_{k=1}^p | u_k | \bigg) \leq C \Vert h \Vert_{m, \alpha; **} \]
\end{Lemma}
\begin{proof}
Assume again $n>3$.
Only the second inequality has to be shown.
First observe that since $|(h-u_k)(c_k)|$ is minimal, we have $(h-u_k)(c_k) \perp u_k(c_k)$ and hence
\[  |(u_k) (c_k)| \leq |h (c_k)| \leq C \Vert h \Vert_{m, \alpha} \leq C \Vert h \Vert_{m, \alpha; **}. \]
It remains to bound $\Vert \bar h \Vert_{m, \alpha; *}$.
Let $h = \bar h' + \sum_{k=1}^p \rho_k u_k'$ be an arbitrary decomposition of $h$ analogous to ($*$).
We will show that
\[ \Vert \bar h \Vert_{m, \alpha; *} \leq C \bigg( \Vert \bar h' \Vert_{m, \alpha; *} + \sum_{k=1}^p |u_k' | \bigg). \]
Note that by the minimal choice of $u_k$, we have $|(u_k - u_k') (c_k)| \leq | (h - u_k') (c_k) |$.
We now use
\[ \Vert \bar h \Vert_{m, \alpha; *} \leq \Vert \bar h' \Vert_{m, \alpha; *} + \sum_{k=1}^p \Vert \rho_k(u_k - u'_k) \Vert_{m, \alpha; *}  \]
and bound the last term by $C \sum_{k=1}^p M_k |(u_k - u'_k)(c_k)| \leq C \sum_{k=1}^p M_k |(h-u'_k)(c_k)|$ where 
\[ M_k = \max_{\widehat N_k} W^{-1} = W^{-1}(c_k). \]
So if $\rho_k(c_k)=1$, then $M_k |(h-u'_k) (c_k)| \leq C \Vert h- \sum_{l=1}^p \rho_{l} u'_{l} \Vert_{m,\alpha;*} = \Vert \bar h' \Vert_{m, \alpha;*}$.
If not, we have a uniform bound on $R_k$, hence on $M_k$ and $M_k |(h-u'_k)(c_k)| \leq C \Vert \bar h' \Vert_{m, \alpha}$.
\end{proof}

\section{Application of the inverse function theorem} \label{sec:Applift}
We will use the following estimate on $L^{-1}$ which we will prove in the next section:
\begin{Proposition} \label{Prop:Linv}
There are $R_0 = R_0(n, V), \Lambda = \Lambda(n, V) < \infty$ such that whenever $\vol M_{hyp} < V$ and $R_{\min} > R_0$, then the operator $L = L_{g_{\overline{\sigma}}} : C^{m, \alpha}(M_{\overline{\sigma}}; \Sym_2 T^*) \to C^{m-2, \alpha}(M_{\overline{\sigma}}; \Sym_2 T^*)$ is invertible and
\[ \Vert h \Vert_{m, \alpha; **} \leq \Lambda \Vert L_{g_{\overline{\sigma}}} h\Vert_{m-2, \alpha; *} \]
for any symmetric bilinear form $h$ on $M_{\overline{\sigma}}$.
\end{Proposition}
Observe that there are different types of norms on both sides of this inequality.
Thus in order to construct a perturbation of $g_{\overline{\sigma}}$ which is Einstein we cannot simply use this estimate to strictly apply the inverse function theorem on Banach spaces.
However, we will show that the trivial Einstein deformations which make the difference between these two norms, have a weak influence on the nonlinear term of the equation we want to solve.

We will now prove Theorem \ref{Thm:main} assuming Proposition \ref{Prop:Linv}.
\begin{proof}[Proof of Theorem \ref{Thm:main}]
We only consider the case $n>3$.
It will be clear how to adapt the proof to the $3$ dimensional case by considering the $N'_k$ as well.

In the following we set $M = M_{\overline{\sigma}}$, $g = g_{\overline{\sigma}}$, $\Phi = \Phi_g$ and $L = L_g$.

Assume that $R_{\min} > R_0$.
We want to find $h \in C^{m, \alpha}(M; \Sym_2 T^*)$ such that the equation $\Phi(g + h) = 0$ holds.
It will then follow from elliptic regularity that $h$ is actually smooth.
The equation is equivalent to the fixed point equation 
\[ h = \Psi(h) = h - L^{-1} \Phi( g + h ). \] 
In order to solve this equation for large $R_{\min}$, it suffices to show that there is an $\varepsilon = \varepsilon(n, V) > 0$ such that $\Psi$ is $\frac12$-Lipschitz with respect to the $\Vert \cdot \Vert_{m, \alpha; **}$-norm on $B_\varepsilon = \{ h \in C^{m, \alpha}(M; \Sym_2 T^*) \;\; : \;\; \Vert h \Vert_{m, \alpha; **} < \varepsilon \}$.
Then, assuming $R_{\min}$ to be large enough, we can achieve $\Vert \Psi(0) \Vert_{m, \alpha; **} = \Vert L^{-1} \Phi(g) \Vert_{m, \alpha; **} \leq \Lambda \Vert \Phi(g) \Vert_{m-2, \alpha; *} < \frac12 \varepsilon$ and apply Banach's fixed point theorem.

For $h_0, h_1 \in B_\varepsilon$ and $h_t = (1-t) h_0 + t h_1$ we compute
\begin{multline*} \left\Vert \Psi(h_0) - \Psi(h_1) \right\Vert_{m, \alpha; **} 
\leq \left\Vert \int_0^1 L^{-1} \left( L - \dd\Phi_{g+h_t} \right)(h_0 - h_1) \dd t \right\Vert_{m, \alpha; **} \\
\leq \Lambda \int_0^1 \left\Vert \left( \dd\Phi_g - \dd\Phi_{g+h_t} \right) (h_0 - h_1) \right\Vert_{m-2, \alpha; *} \dd t
\end{multline*}
(Observe that the subscript of $\dd\Phi$ now indicates the point at which the derivative is taken rather than the background metric which we used to define $\Phi$.)
Thus, it suffices to show that for any $h \in B_\varepsilon$ and $h' \in C^{m, \alpha}(M; \Sym_2 T^*)$ we have
\[ \left\Vert \dd\Phi_g(h') - \dd\Phi_{g+h} (h') \right\Vert_{m-2, \alpha; *} \leq \delta(\varepsilon) \Vert h' \Vert_{m, \alpha; **} \]
for some universal $\delta(\epsilon)$ with $\delta(\varepsilon) \to 0$ as $\varepsilon \to 0$.

Represent $h = \bar h + \sum_{k=1}^p \rho_k u_k$ and $h' = \bar h' + \sum_{k=1}^p \rho_k u'_k$ where the $u_k, u'_k$ are trivial Einstein variations.
Then
\begin{multline*}
 \left\Vert \dd\Phi_g(h') - \dd\Phi_{g+h} (h') \right\Vert_{m-2, \alpha; *}  \\
\leq \left\Vert \dd\Phi_g(\bar h') - \dd\Phi_{g+h} (\bar h') \right\Vert_{m-2, \alpha; *} + \sum_{k=1}^p \left\Vert \dd\Phi_g(\rho_k u'_k) - \dd\Phi_{g+h} (\rho_k u'_k) \right\Vert_{m-2, \alpha; *}.
\end{multline*}
The first term can immediately be bounded by $C \Vert h \Vert_{m, \alpha} \Vert \bar h' \Vert_{m, \alpha; *} \leq C' \Vert h \Vert_{m, \alpha; **} \Vert \bar h' \Vert_{m, \alpha; *}$.
As for the second term we have
\begin{multline*} \left\Vert \dd\Phi_g(\rho_k u'_k) - \dd\Phi_{g+h} (\rho_k u'_k) \right\Vert_{m-2, \alpha; *} \\
\leq \left\Vert \dd\Phi_g(\rho_k u'_k) - \dd\Phi_{g+\rho_k u_k} (\rho_k u'_k) \right\Vert_{m-2, \alpha; *} 
+ \left\Vert \dd\Phi_{g + \rho_k u_k} (\rho_k u'_k) - \dd\Phi_{g+\bar h + \rho_k u_k} (\rho_k u'_k) \right\Vert_{m-2, \alpha; *}.
\end{multline*}
Now, since $u_k$ is a trivial Einstein variation, we can use (\ref{eq:trivEinsnoinfl}) to bound the first term by $C | u_k | | u'_k |$.
The second term is bounded by $C \Vert \bar h \Vert_{m, \alpha; *} | u'_k |$.
We conclude
\begin{multline*}
 \left\Vert \dd\Phi_g(h') - \dd\Phi_{g+h} (h') \right\Vert_{m-2, \alpha; *}
 \leq C \Vert h \Vert_{m, \alpha; **} \Vert \bar h' \Vert_{m, \alpha; *} + C \sum_{k=1}^p \left( |u_k| |u'_k| + \Vert \bar h \Vert_{m,\alpha; *} |u'_k| \right) \\
\leq  C \Vert h \Vert_{m, \alpha; **} \Vert \bar h' \Vert_{m, \alpha; *} + C \big( \Vert \bar h \Vert_{m, \alpha; *} + \sum_k |u_k| \big) \sum_l |u'_l|.
\end{multline*}
By an appropriate choice of $\bar h$ and $u_k$, the right hand side can be made arbitrarily close to $C \Vert h \Vert_{m, \alpha; **} \left( \Vert \bar h' \Vert_{m, \alpha; *} + \sum_l |u'_l| \right)$ what in turn by a good choice of $\bar h'$ and $u'_k$ can be made arbitrarily close to $C \Vert h \Vert_{m, \alpha; **} \Vert h' \Vert_{m, \alpha; **} \leq C \varepsilon \Vert h' \Vert_{m, \alpha; **}$.
This proves the desired bound and hence the theorem.
\end{proof}

\section{Estimates for $L^{-1}$} \label{sec:estforLm}
This section will be occupied with the proof of Proposition \ref{Prop:Linv}.
For the sake of a clear exposition of the main ideas we will defer most of the technical arguments to sections \ref{sec:notoncusp} and \ref{sec:notonBH}.
We first establish a bound on the $\Vert \cdot \Vert_{m, \alpha}$-norm:
\begin{Lemma} \label{Lem:Linvsimple}
There are $R_0 = R_0(n, V), \Lambda = \Lambda(n, V) < \infty$ such that whenever $\vol M_{hyp} < V$ and $R_{\min} > R_0$, then the operator $L_{g_{\overline{\sigma}}} : C^{m, \alpha}(M_{\overline{\sigma}}; \Sym_2 T^*) \to C^{m-2, \alpha}(M_{\overline{\sigma}}; \Sym_2 T^*)$ is invertible and
\[ \Vert h \Vert_{m, \alpha} \leq \Lambda \Vert L_{g_{\overline{\sigma}}} h\Vert_{m-2, \alpha; *} \]
for any symmetric bilinear form $h$ on $M_{\overline{\sigma}}$.
\end{Lemma}
\begin{proof}
The proof of this Lemma is similar to that of \cite[Proposition 3.2]{And}.
Observe that in this Proposition the right hand side of the inequality reads $\Lambda \log R_{\max} \Vert L_{g_{\overline{\sigma}}} h \Vert_{m-2,\alpha}$.
In our case, we can don't need the $\log R_{\max}$ factor, but have to make use of a stronger norm of $L_{g_{\overline{\sigma}}}$.

Recall that we have the Schauder estimate
\[ \Vert h \Vert_{m,\alpha} \leq C( \Vert L_{g_{\overline{\sigma}}} h \Vert_{m-2,\alpha} + \Vert h \Vert_{0} ) \]
where $C$ is uniform.
So it is enough to show that: \vspace{0.5\baselineskip}

\it There are $R_0 = R_0(n, V), \Lambda'=\Lambda'(n, V) < \infty$ such that whenever $\vol M_{hyp} < V$ and $R_{\min} > R_0$, we have
\[ \Vert h \Vert_{0}  \leq  \Lambda' \Vert L_{g_{\overline{\sigma}}} h \Vert_{m-2,\alpha; *} \]
for all symmetric bilinear forms $h$ on any $M_{\overline{\sigma}}$. \rm \vspace{0.5\baselineskip}

Assume that this statement was wrong.
Then we can find a sequence of hyperbolic manifolds $M_{hyp}^i$ with basepoints $y^i \in M_{thick}^i$ and $\vol M^i_{hyp}$ uniformly bounded from above as well as a sequences of $\overline{\sigma}^i$ such that $R_{\min}^i \to \infty$ and symmetric bilinear forms $h^i$ on $M^i = M_{\overline{\sigma}^i}$ such that for $g^i = g_{\overline{\sigma}^i}$, $L^i = L_{g^i}$ and $f^i := L^i h^i$
\[ \Vert h^i \Vert_{0} = 1, \qquad \text{but} \qquad \Vert f^i \Vert_{m-2,\alpha; *} \longrightarrow 0 \]
 as $i \to \infty$.
So there are points $x^i \in M^i$ such that $|h^i| (x^i) > \gamma$ for some universal $\gamma > 0$.
The Schauder estimate gives us a uniform $C^{m, \alpha}$-bound for the $h^i$.\\[0.1\baselineskip]

$1^\circ \quad$ In the first step we show that there are sequences $d^i \to \infty$ and $w^i \to 0$ such that $| h^i | < w^i$ on $B_{d^i} (y^i)$.

Consider an arbitrary subsequence of counterexamples.
After passing to a subsequence again, the pointed Riemannian manifolds $(M^i, y^i)$ Gromov-Hausdorff converge to a pointed hyperbolic manifold $(M^\infty_{hyp}, y^\infty)$ of finite volume.
Furthermore, the $h^i$ subconverge to a symmetric bilinear form $h^\infty$ on $M_{hyp}^\infty$ such that $L^\infty (h^\infty) = 0$ (here $L^\infty = L_{g_{hyp}^\infty}$).

Denote by $(M_{hyp}^\infty)^s$ the manifold obtained from $M_{hyp}^\infty$ by truncating its cusps at distance $s$ from the basepoint $y^\infty$.
Using Stoke's theorem and (\ref{eq:hypWbck}), we find
\[ \int_{(M_{hyp}^\infty)^s} |dh^\infty|^2 + |\DIV h^\infty|^2 + (n-2) |h^\infty|^2 + (\tr h^\infty)^2 = \int_{\partial (M_{hyp}^\infty)^s} Q(h^\infty, \nabla h^\infty) \]
where the right hand side goes to $0$ as $s \to \infty$.
So $h^\infty \equiv 0$ and we conclude that for any $d$ we have $|h^i| \to 0$ uniformly on $B_d(y^i)$ for a subsequence.
Since we started with an arbitrary subsequence, this implies that for any $d$ we have $|h^i| \to 0$ on $B_d(y^i)$ uniformly for the \emph{whole} sequence and hence the claim.\\[0.1\baselineskip]

$2^\circ \quad$ Next, we give an estimate for $h^i$ on the $\widehat N_{k}^i$ (and ${N'}_{k}^i$ in dimension $3$).

Choose coordinates $(r, \theta, x_3, \ldots, x_n)$ on these components (to be precise on their universal covers).
Observe that $\widehat N_{k}^i \setminus B_1(T_{k}^i)$ carries the exact black-hole metric.
We have $|f^i| < \Vert f^i \Vert_{m-2, \alpha; *} W$.
Since $\dist(y^i, T_{k}^i) < \diam M_{thick}^i$ is uniformly bounded, we find that $| h^i | < w^i \to 0$ around the boundaries of the $\widehat N_{k}^i$.

Consider the restriction of $h^i$ and $f^i$ to some $\widehat N_{k}^i$ and take their average under the $S^1 \times \IR^{n-2}$-action, i.e. let $T^{n-1}(r') := \widehat N_{k}^i(r=r')$ be the cross-sectional torus at the coordinate $r = r'$ and set
\[ \widehat h_{st}^i(r) := \frac1{\vol T^{n-1}(r)} \int_{T^{n-1}(r)} h_{st}^i. \]
Analogously define $\widehat f^i$.
Obviously, $\widehat h^i$ and $\widehat f^i$ are $S^1 \times \IR^{n-2}$ invariant and $L^i \widehat h^i = \widehat f^i$.
Furthermore, still 
\[ |\widehat f^i| < \Vert f^i \Vert_{m-2, \alpha; *} W = \Vert f^i \Vert_{m-2, \alpha; *} \bigg[ \Big( \frac{r}{R_{k}^i}  \Big)^{0.1} + r^{-0.1}\bigg] \]
and since $\nabla h^i$ is uniformly bounded and $\diam T^{n-1} (r) < C \frac{r}{R_{k}^i}$, we conclude
\[ | \widehat h^i - h^i | < C \frac{r}{R_{k}^i} \qquad \text{on} \quad T^{n-1}(r). \]
We can now apply Proposition \ref{Prop:uglyestimate} to conclude
\begin{equation} | h^i | < C \left( w^i + \Vert f^i \Vert_{m-2, \alpha; *} + r^{-n+1.1} + \frac{r}{R_{k}^i} \right). \label{eq:estimateonnose} \end{equation}\\[0.1\baselineskip]

$3^\circ \quad$ We can make the following conclusions on $x^i$:
From $1^\circ$ we already know that $\dist(y^i, x^i) \to \infty$.
This implies that $x^i$ eventually lies in some $\widehat N_{k}^i$ (or ${N'}_{k}^i$ in dimension $3$) and $\frac{r(x^i)}{R_{k}^i} \to 0$.
So by (\ref{eq:estimateonnose}) we conclude that $r(x^i)$ has to stay bounded.
This means that the $x^i$ have to stay in bounded distance to some core tori $(\widehat T^{n-2}_{k})^i$ of $\widehat N_{k}^i$ (or of ${N'}_{k}^i$).

So there is a sequence ${d'}^i$ such that the universal covers $(\widetilde{B}_{{d'}^i}(x^i), x^i)$ Gromov-Hausdorff subconverge to the black-hole metric $(M_{BH}, x^\infty)$ and the $h^i$ subconverge to some $h^\infty$ on $M_{BH}$ which satisfies $L^\infty h^\infty = 0$ and $h^\infty(x^\infty) \not= 0$.
Moreover, since the the pointed manifolds $(\widehat N_{k}^i, x^i)$ collapse to a ray, $h^\infty$ is invariant under the $S^1 \times \IR^{n-2}$-action.
From (\ref{eq:estimateonnose}) we also conclude that $| h^\infty | < C r^{-n+1.1}$.

We can now use Proposition \ref{Prop:EinstvarBH} to find that $h^\infty \equiv 0$, a contradiction.
\end{proof}

Finally, we can use Lemma \ref{Lem:Linvsimple} to refine our result and prove Proposition \ref{Prop:Linv}:
\begin{proof}[Proof of Proposition \ref{Prop:Linv}]
Analogously to the proof of Lemma \ref{Lem:Linvsimple}, we assume that the hypothesis was wrong and that we have sequences $M^i, \overline{\sigma}^i, h^i$ such that $R_{\min}^i \to \infty$ and
\[ \Vert h^i \Vert_{0; **} = 1, \qquad \text{but} \qquad \Vert f^i \Vert_{m-2,\alpha; *} \longrightarrow 0 \]
for $f^i = L^i h^i$ (we also used (\ref{eq:ssSchauder}) here).

By Lemma \ref{Lem:Linvsimple} we have $\Vert h^i \Vert_{0} \to 0$.
We now change the $h^i$ by certain trivial Einstein variations of the $\widehat N_{k}^i$ (or ${N'}_{k}^i$ in dimension $3$):
Let $u_{k}^i$ be those trivial Einstein deformations as obtained in Lemma \ref{Lem:trivcenter} and set $\bar h^i = h^i - \sum_{k=1}^{p^i} \rho_{k}^i u_{k}^i$.
Then, since $\Vert h^i \Vert_0 \to 0$ we have $|u_{k}^i| \to 0$ as $i \to \infty$ and by Lemma \ref{Lem:trivcenter}
\[ \Vert \bar h^i \Vert_{**; 0} \leq \Vert \bar h^i \Vert_{0; *} \leq C \Vert \bar h^i \Vert_{0; **} \]
for some uniform constant $C$.
So we conclude that $\frac12 < \Vert \bar h^i \Vert_{0; **} < \frac32$ for large $i$ and hence we have the uniform estimate $c < \Vert \bar h^i \Vert_{0; *} < C$.
However, we still have $\Vert \bar h^i \Vert_0 \to 0$.
Finally, setting $\bar f^i = L^i \bar h^i$, we get $\Vert \bar f^i \Vert_{m-2, \alpha; *} \to 0$.

By the lower bound on $\Vert \bar h^i \Vert_{0; *}$, we can find points $x^i \in M^i$ such that
\[ W^{-1}(x^i) |\bar h^i|(x^i) > \gamma > 0. \]
Since $\Vert \bar h^i \Vert_0 \to 0$, we conclude $W(x^i) \to 0$.
So the $x^i$ eventually lie in certain $\widehat N_{k^i}^i$ (or ${N'}_{k^i}^i$), $R_{k^i}^i \to \infty$ and the distance of the $x^i$ to both $T_{k^i}^i$ as well as $(\widehat T_{k^i}^{n-2})^i$ goes to infinity.
So there is a sequence ${d'}^i$ such that the universal covers $(\widetilde{B}_{{d'}^i}(x^i), x^i)$ converge to hyperbolic space $(\IH^n, x^\infty)$ on which we can choose coordinates $(r^\infty, x_2, \ldots, x_n)$ with $r^\infty(x^\infty) = 1$ and $\frac{r}{r^i} \to r^\infty$ where $r^i := r(x^i)$ (observe that we choose those coordinates in which the hyperbolic metric takes the form (\ref{eq:hypmetricrep})).
In order to analyze the limiting behaviour of $\bar h^i$, we have to distinguish three cases:\\[0.1\baselineskip]

$1^\circ \quad$ For a subsequence we have $r^i (R_{k^i}^i)^{-1/2} \to \infty$.

Then we have the (local) convergence
\[ \left( \frac{R_{k^i}^i}{r^i} \right)^{0.1} W = \left( \frac{r}{r^i} \right)^{0.1} + \left(\frac{ R_{k^i}^i}{(r^i)^2} \cdot \frac{r^i}{r} \right)^{0.1} \longrightarrow (r^\infty)^{0.1}. \]
So $( \frac{R_{k^i}^i}{r^i} )^{0.1}  \bar h^i$ is locally bounded and $( \frac{R_{k^i}^i}{r^i} )^{0.1} \bar f^i \to 0$ locally.
Hence the $\bar h^i$ subconverge to some nonzero $\bar h^\infty$ on $\IH^n$ which satisfies $|\bar h^\infty| < C (r^\infty)^{0.1}$ and $L^\infty \bar h^\infty = 0$ (observe here that by Schauder estimates, we have uniform local bounds on some derivatives of the $\bar h^i$).
Since the sequence $(B_{{d'}^i}(x^i), x^i)$ collapses to a line, $\bar h^\infty$ must be invariant under the group $\IR^{n-1}$ acting on the last coordinates.
We can now use Proposition \ref{Prop:notonhyp} to obtain a contradiction.\\[0.1\baselineskip]

$2^\circ \quad$ For a subsequence we have $r^i (R_{k^i}^i)^{-1/2} \to 0$.

This time we have the convergence
\[ (r^i)^{0.1} W = \left( \frac{(r^i)^2}{R_{k^i}^i} \cdot \frac{r}{r^i} \right)^{0.1} + \left( \frac{r}{r^i} \right)^{-0.1} \longrightarrow (r^\infty)^{-0.1}. \]
Now we can use the same arguments as in $1^\circ$ to construct $\bar h^\infty$ on $\IH^n$ which obeys the bound $|\bar h^\infty| < C (r^\infty)^{-0.1}$.
This also contradicts Proposition \ref{Prop:notonhyp}.
\\[0.1\baselineskip]

$3^\circ \quad$ For a subsequence we have $r^i (R_{k^i}^i)^{-1/2} \to q$ where $0 < q < \infty$.

This means that the points $x^i$ stay within bounded distance to the $c_{k^i}^i$.
Let $c^\infty \in \IH^n$ be one of their limit points.
We have the convergence
\[ (r^i)^{0.1} W =  \left( \frac{(r^i)^2}{R_{k^i}^i} \cdot \frac{r}{r^i} \right)^{0.1} + \left( \frac{r}{r^i} \right)^{-0.1} \longrightarrow q^{0.2} (r^\infty)^{0.1} + (r^\infty)^{-0.1}. \]
Hence the same reasoning as in $1^\circ$ yields a nonzero $\bar h^\infty$ which satisfies $| \bar h^\infty | < C( (r^\infty)^{0.1} + (r^\infty)^{-0.1} )$.
So by Proposition \ref{Prop:notonhyp}, $\bar h^\infty$ must be trivial.

However, by the construction of the $\bar h^i$ we get that  $|(\bar h^\infty)(c^\infty) | \leq | (\bar h^\infty - u)(c^\infty) |$ for any trivial Einstein variation $u$, contradicting the fact that $\bar h^\infty$ is nonzero.
\end{proof}

\section{Einstein variations of the hyperbolic cusp metric} \label{sec:notoncusp}
Consider the hyperbolic metric
\[ g_{hyp} = r^{-2} \dd r^2 + r^2 (\dd x_2^2 + \ldots + \dd x_n^2) \]
on $\IR^+ \times \IR^{n-1}$ and the parabolic isometric action of $\IR^{n-1}$ by translations on the second factor.

Set $L := L_{g_{hyp}}$.
We will prove the following result:
\begin{Proposition} \label{Prop:notonhyp}
Let $h$ be a symmetric bilinear form on $\IH^n$ that is invariant under the $\IR^{n-1}$-action.
Assume furthermore that $|h| < r^{0.1} + r^{-0.1}$. \\
Then $L h = 0$ implies that $h$ is trivial.

Thus, if even $|h| < r^{\pm 0.1}$, then $h \equiv 0$.
\end{Proposition}
\begin{proof}
We assume $|h| < r^{0.1} + r^{-0.1}$.
Express $h = h_{ij} \dd x_i \dd x_j$ where we set $x_1 = r$.
Then the $h_{ij}$ only depend on $r$ and the bound on $|h|$ implies
\[ r^2 |h_{11}|(r), \quad |h_{1i}|(r), \quad r^{-2} |h_{ij}|(r) < r^{0.1} + r^{-0.1} \]
for $i, j > 1$.

The equation $Lh = 0$ writes out as (see (\ref{eq:hypWbck}))
\[  \triangle h + 2 h - 2 (\tr h) g_{hyp} = 0 \]
which implies
\begin{alignat}{4}
r^2 (r^2 h_{11})&'' & + n r (r^2 h_{11} )&'& - 2(n-1) (r^2 &h_{11}) &&= 0 \tag{I} \\
r^2 h_{1i}&'' & + n r h_{1i}&'& - n &h_{1i} &&= 0 \tag{II} \\
r^2 (r^{-2} h_{ij} )&'' & + n r ( r^{-2} h_{ij} )&'& - 2 \delta_{ij} \sum_{k=2}^n r^{-2} &h_{kk} &&= 0 \tag{III}
\end{alignat}
The trace of $h$ satisfies (see (\ref{eq:tr}))
\[  \triangle \tr h - 2(n-1)\tr h = 0. \]
In terms of coordinates, this implies for $q(r) = \tr h(r)$
\begin{equation}  r^2 q'' + n r q' - 2(n-1) q = 0. \tag{IV} \end{equation}
The solutions of (I) and (IV) are both of the form $A_1 r^{\gamma_1} + A_2 r^{\gamma_2}$ with $\gamma_{1/2} = \frac12(-n+1 \pm \sqrt{n^2 + 6n - 7})$.
Hence by the bound on $|h|$ we get $r^2 h_{11} \equiv \tr h \equiv 0$ and plugging this into (III) gives 
\[ r^2 (r^{-2} h_{ij})'' + n r (r^{-2} h_{ij})' = 0. \]
Solutions of this equation are of the form $h_{ij}(r) = A_1 r^2 + A_2 r^{-n+3}$ and thus $h_{ij} = u_{ij} r^2$ (for $i,j > 1$).

Finally, (II) implies $h_{1i}(r) = A_1 r + A_2 r^{-n}$, hence $h_{1i} \equiv 0$ (for $i > 1$).
\end{proof}

\section{Variations of the black-hole metric} \label{sec:notonBH}
Consider the black-hole metric
\[ g = g_{BH} = V^{-1} \dd r^2 + V \dd \theta^2 + r^2 (\dd x_3^2 + \ldots + \dd x_n^2) \]
on $M_{BH} \approx \IR^2 \times \IR^{n-2}$.
Set $L = L_{g}$.
Recall that $g$ is asymptotic to the hyperbolic metric
\[ g_{hyp} = r^{-2} \dd r^2 + r^2 (\dd \theta^2 + \dd x_2^2 + \ldots + \dd x_n^2) \]
for $r \to \infty$ in the sense that $| \nabla^m (g - g_{hyp}) | = O( r^{-n+1} )$.
This is why we can estimate
\[ | L_{g} h - L_{g_{hyp}} h | \lesssim O (r^{-n+1}) | h | + O ( r^{-n+1} ) | \nabla h | + O( r^{-n+1} ) | \nabla^2 h | \]
for $r \to \infty$.

In the following we will analyze Einstein variations of $g_{BH}$ or variations which are almost Einstein.
We will hereby always assume that these variations are invariant under the $S^1 \times \IR^{n-2}$ action.
When we compare $g_{BH}$ with $g_{hyp}$, this action becomes the parabolic $\IR^{n-1}$ action.

We remark that Olivier Biquard has independently found elementary proofs of some of the following results (\cite{Biq}).

\begin{Proposition} \label{Prop:uglyestimate}
Let $R > r_+$ and assume that on $M_{BH}(r \leq R)$ we have $L h = f$ for $S^1 \times \IR^{n-2}$ invariant $h$ and $f$ satisfying $| h |(r) < 1$ and 
\[ | f |(r) < \alpha \left[ \left( \frac{r}{R} \right)^{0.1} + r^{-0.1} \right] \]
for all $r \leq R$ and some $\alpha < 1$.
Then
\[ |h|(r) < C \left( |h| (R) + \alpha  + r^{-n+1.1} \right) \]
for some universal constant $C$ (which is independent of $R$).
\end{Proposition}

We will need a technical Lemma.
Note that from now on whenever we use the notation $O( \varphi(r) )$ for a function $\varphi(r)$, we indicate an error term whose absolute value is \emph{always} (not only for $r \to \infty$) bounded above by $C \varphi(r)$ where $C$ is a universal constant.
\begin{Lemma} \label{Lem:ODE2est}
Let $a, b \in \IR$ and $0 \leq B_1, B_2 \leq \infty$.
Consider a solution $f : (B_1, B_2) \to \IR$ of the ODE 
\[ r^2 f''(r) + a r f'(r) + b f(r) = \varphi(r) \]
for some $\varphi : (B_1, B_2) \to \IR$.
Assume that $a, b$ are chosen in such a way that the corresponding homogeneous ODE (for $\varphi \equiv 0$) has the general solution $f(r) = A_1 r^{\gamma_1} + A_2 r^{\gamma_2}$ with $\gamma_1, \gamma_2 \in \IR$ and $\gamma_1 < \gamma_2$.

Now, suppose $\varphi(r) = \sum_{k=1}^p O( r^{\delta_k} )$ where we assume that $\delta_k \not= \gamma_1, \gamma_2$ for each $k$.
Then $f(r) = A_1 r^{\gamma_1} + A_2 r^{\gamma_2} + \sum_{k=1}^p O( r^{\delta_k} )$.

Here the coefficients in $O( r^{\delta_k} )$ only depend on $a, b, \delta$ and the coefficients in the error terms of $\varphi$.
\end{Lemma}
\begin{proof}
The Lemma follows by simple integration.
\end{proof}

\begin{proof}[Proof of Proposition \ref{Prop:uglyestimate}]
We assume from now on that $|h|(r) < 1$ and $R > r_+ + 2$.
Using the Schauder estimates we find that this implies $|\nabla^l h | < C_l$, so
\begin{equation} | L_{g} h - L_{g_{hyp}} h | = O( r^{-n + 1} ) \label{eq:Lhypest} \end{equation}
for $r > r_+ + 1$.
In coordinates, the bound on $h$ implies that 
\[ r^2 |h_{11}|(r), \quad |h_{1i}|(r), \quad r^{-2} |h_{ij}|(r) < C \]
where $i,j > 1$ and $r > r_+ + 1$ for the first quantity and $r \geq r_+$ for the rest.

We will use the equations from the last section to derive a better estimate on $h$.
Set $H = |h|(R)$. 

We first show how to bound $r^2 h_{11}$.
By equation (I) of the last section and (\ref{eq:Lhypest}) it satisfies
\[ r^2 (r^2 h_{11})'' + n r (r^2 h_{11})' - 2 (n-1) (r^2 h_{11}) = r^2 f_{11} + O( r^{-n+1} ), \]
where $r^2 f_{11}(r) = O(\alpha ( \frac{r}{R} )^{0.1} ) + O (\alpha r^{-0.1})$ for $r > r_+ + 1$.
Lemma \ref{Lem:ODE2est} gives us
\[ r^2 h_{11} (r) = A_1 r^{\gamma_1} + A_2 r^{\gamma_2} +  O ( \alpha ( \tfrac{r}{R} )^{0.1} ) + O ( \alpha r^{-0.1} ) + O( r^{-n+1} ), \]
where $\gamma_{1/2} = \frac12 (-n+1 \pm \sqrt{n^2+6n-7})$.
Observe that $\gamma_1 > 0.1$ and $\gamma_2 < -n+1$.
Since $r^2 h_{11}$ and the error terms above are bounded for say $r \in (r_+ +1, r_+ + 2)$, we conclude that $|A_2| < C$ for some universal $C$.
For $r = R$, we furthermore obtain $|A_1| < C H R^{-\gamma_1} + O (R^{\gamma_2 - \gamma_1}) + O(\alpha R^{-\gamma_1}) + O(R^{-n+1-\gamma_1})$.
Thus
\begin{equation} \label{eq:h11bound} 
r^2 |h_{11}|(r) < C \bigg[ H \Big( \frac{r}{R} \Big)^{\gamma_1} + \alpha \Big( \frac{r}{R} \Big)^{0.1} + \alpha r^{-0.1} + r^{-n+1} \bigg] < C (H + \alpha + r^{-n+1}). 
\end{equation}
Using (IV) from the last section, we conclude that the same bound holds for $\tr h$.
Moreover, we can estimate $h_{1i}$ for $i>1$ by the same method (this time we have to use (II) and the fundamental solutions are $r$ and $r^{-n}$).

Using the first estimate from (\ref{eq:h11bound}), we can now bound $r^{-2} h_{ij}$ for $i, j > 1$.
By (III) we obtain for $r > r_+ + 1$
\[  r^2 (r^{-2} h_{ij})'' + n r (r^{-2} h_{ij})' 
 =    O({ \textstyle H ( \frac{r}{R} )^{\gamma_1}}) + O ( \alpha ({ \textstyle \frac{r}{R} } )^{0.1} ) + O(\alpha r^{-0.1}) + O( r^{-n+1} )  . \]
Thus using $O(r^{-n+1}) < O(r^{-n+1.1})$, we conclude from Lemma \ref{Lem:ODE2est}
\[ r^{-2} h_{ij}(r) = A_1 + A_2 r^{-n+1} \\
 + O({ \textstyle H ( \frac{r}{R} )^{\gamma_1}}) + O ( \alpha ({ \textstyle \frac{r}{R} } )^{0.1} ) + O(\alpha r^{-0.1}) + O( r^{-n+1.1} ) . \]
As before, we find that $|A_2| < C$ and setting $r=R$ yields $|A_1| < C H + O(H + \alpha) + O(R^{-n+1.1})$, so
\[ r^{-2} | h_{ij} | (r) < C \left( H + \alpha + r^{-n+1.1} \right). \qedhere \] 
\end{proof}

We will now prove the second result of this section.
\begin{Proposition} \label{Prop:EinstvarBH}
Let $h$ be an $S^1 \times \IR^{n-2}$ invariant Einstein variation of $g_{BH}$ and assume $|h|(r) \to 0$ for $r \to \infty$.
Then $h \equiv 0$.
\end{Proposition}
We note that with a little more work, it is even possible to deduce that any $S^1 \times \IR^{n-2}$ invariant Einstein variation $h$ with $|h|(r) \to 0$ as $r \to \infty$ is of the form
\[ h = - \tr u\frac{n-1 }{V r^{n-1}} \dd r^2 - \tr u \frac{V V'}{2 r} \dd \theta^2 + 2 (\tr u) r^{-n+3} (\dd x_3^2 + \ldots + \dd x_n^2 ) + u_{ij} r^2 \dd x_i \dd x_j \]
for some symmetric $(n-2) \times (n-2)$ matrix $u_{ij}$ indexed by $i,j = 3, \ldots, n$.

Assume from now on that $|h|(r) \to 0$ as $r \to \infty$ and that $L h = 0$.
Using Proposition \ref{Prop:uglyestimate}, we find that we even have $|h|(r) < C r^{-n+1.1}$.
By Schauder's estimates we can deduce the same decay for all covariant derivatives of $h$.

\begin{Lemma} \label{Lem:trDIVis0}
We have $\tr h \equiv 0$, $\DIV h \equiv 0$ and hence $\dd \Ric_{g}(h) + (n-1)h = 0$.
\end{Lemma}
\begin{proof}
This follows from the maximum principle applied to (\ref{eq:tr}) resp. (\ref{eq:bg}) and the fact that $\tr h$ and $\beta( h)$ are decaying.
\end{proof}

\begin{Lemma}
 We have $h_{1i} = h_{i1} \equiv 0$ for all $i \geq 2$.
\end{Lemma}
\begin{proof}
 Writing out the equation $\DIV h = 0$ in terms of the $h_{ij}$ gives for $i \geq 2$
\[ 0 = -(\DIV h)_i = V h'_{1i} + \Big( V' + (n-2) \frac{V}r \Big) h_{1i} \]
The solutions of these ODEs behave like $\frac1{r-r_+}$ for $r \to r_+$, so the $h_{1i}$ must be constantly zero.
\end{proof}

Now we will alter $h$ by an infinitesimal diffeomorphism $\DIV^*_g \xi$ for some $1$-form $\xi$ to eliminate its $11$ entry.
Observe that by (\ref{eq:Ricdiffeo}) for every 1-form $\xi$ we have
\[ \dd \Ric_{g}(\DIV^*_{g} \xi) + (n-1) \DIV^*_{g} \xi = 0 \]
since $g$ is Einstein.
So for any $1$-form $\xi$ the bilinear form $h + \DIV^*_{g} \xi$ will still be an infinitesimal Einstein variation.
However, we might lose the divergence or trace freeness.
\begin{Lemma} \label{Lem:xi}
There is an $S^1 \times \IR^{n-2}$ invariant 1-form $\xi = \xi_1 (r) dr$ such that for $k = h + \DIV^* \xi$ we have $k_{1i} = k_{i1} = 0$ for $i = 1, \ldots, n$ and $|k |(r) < C$.

Moreover, if $k \equiv 0$, then $\xi \equiv 0$ and hence $k \equiv h$.
\end{Lemma}
\begin{proof}
We compute
\begin{alignat*}{1}
(\DIV^*\xi)_{11} &= \xi'_1 + \frac{V'}{2V} \xi_1, \\ 
(\DIV^*\xi)_{22} &= \tfrac12 V V' \xi_1 \\ 
(\DIV^*\xi)_{ii} &=  r V \xi_1 \qquad \text{for $i \geq 3$} 
\end{alignat*}
The remaining components are zero.

We now solve the ODE $(\DIV^* \xi)_{11} = - h_{11}$.
Observe that it is equivalent to $(V^{1/2} \xi_1)' = - V^{1/2} h_{11}$ and that $V |h_{11}| < C r^{-n+1.1}$.
Hence, the solution
\[ \xi_1 (r) := - \frac1{V^{1/2}} \int_{r_+}^r \frac{V h_{11}}{V^{1/2}} \]
satisfies $V^{1/2}(r) |\xi_1|(r) \leq C (r - r_+)^{1/2}$ which implies smoothness of $\xi$ and boundedness of $\DIV^* \xi$.

Now, if $k \equiv 0$, then $\tr h \equiv 0$ implies $\tr \DIV^* \xi \equiv 0$ and thus
\[ V \xi_1' + \Big( V' + (n-2) \frac{V}r \Big) \xi_1 = 0. \]
Hence $\xi_1(r) = C V^{-1} (r) r^{-n+2}$ which behaves like $\frac{C}{2(n-1)} \frac1{r-r_+}$ as $r \to r_+$ contradicting the smoothness of $\xi$.
\end{proof}
We will now show that $k$ has a very simple form.
In order to do this, we introduce a new coordinate $s = s(r)$ (the distance to the origin) with the property that $s(r_+) = 0$ and $g_{BH} = \dd s^2 + V(r(s)) \dd \theta^2 + r^2(s) (\dd x_3^2 + \ldots + \dd x_n^2)$.
From now on we will only work in the coordinate system $(s, \theta, x_3, \ldots, x_n)$.
Consider a metric $\widehat g$ of the form 
\[ \widehat g(s) = \left( \begin{matrix} 1 & 0 \\ 0 & M (s) \end{matrix}\right). \]
The condition of being Einstein with $\Ric_{\widehat g} = -(n-1) \widehat g$ is equivalent to the following system of ODEs (see e.g. \cite{Lin}):
\begin{alignat}{1}
 \left( \sqrt{\det M} M' M^{-1} \right)' - 2(n-1) \sqrt{\det M} E_{n-1} &= 0 \tag{I} \\
 \chi_{n-2}(M' M^{-1}) - 2(n-2)(n-1) &= 0 \tag{II}
\end{alignat}
where $E_{n-1}$ denotes the unit matrix of rank $n-1$ and $\chi_{n-2}$ the $(n-2)$-th coefficient of the characteristic polynomial, i.e. the elementary symmetric polynomial of degree $2$ in the eigenvalues.
The prime denotes differentiation by $s$.

Now denote by $M = \diag (V(r(s)),r^2(s), \ldots, r^2(s))$ the matrix corresponding to the black-hole metric $g=g_{BH}$ and denote the given Einstein variation corresponding to $k_{2 \leq i,j \leq n}$ by $\dot M = \dot M(s)$.
Then $\dot M$ is a variation $\dot M = \dot M(s)$ of (I) and (II).
We will moreover abbreviate $u = \sqrt{\det M}$ and $\dot u = (\sqrt{\det M})^\cdot = \frac12 u \tr( \dot M M^{-1})$.

\begin{Lemma} \label{Lem:usinh}
 $u = A \sinh (n-1)s$ for some $A > 0$ and $\dot u = \dot A u$.
This implies $\tr_g k = \tr \dot M (M^0)^{-1} \equiv 2 \dot A$.
\end{Lemma}
\begin{proof}
Tracing (I) yields
\[ u'' - (n-1)^2 u = 0 \]
Since $u(0) = 0$, we get $u = A \sinh (n-1)s$ and a variation of this equation gives
$\dot u = \dot A \sinh(n-1)s + \dot B \cosh (n-1)s$.
So $\frac12 \tr k = \dot A A^{-1} + \dot B A^{-1} \cosh (n-1)s / \sinh (n-1)s$.
Since $\tr k$ is bounded, we conclude $\dot B = 0$ and therefore $\tr k \equiv \const$.
\end{proof}

Now observe that by the symmetries $x_i \to - x_i$ ($i \geq 3$), also the matrix
\[ \dot M^\bot := \left( 
\begin{matrix} 
\dot M_{22}   & - \dot M_{23}  & \cdots   & - \dot M_{2n} \\ 
- \dot M_{32} &   \dot M_{33}  & \cdots   &   \dot M_{3n} \\
\vdots        & \vdots         & \ddots   &  \vdots  \\
- \dot M_{n2} &   \dot M_{n3}  & \cdots   &   \dot M_{nn}
\end{matrix}\right) \]
corresponds to an Einstein variation.
Moreover, we can easily see that the entries $\dot M_{23}, \ldots, \dot M_{2n}$ are odd functions in $s$ while all the other entries are even.
So $\dot M$ is invariant under the transformation $\perp$ combined with $s \to -s$.

\begin{Lemma}
 $\dot M = Q M$ where $Q$ is a symmetric matrix with $Q_{2i} = Q_{i2} = 0$ for $i = 2, \dots, n$.
 Hence $\dot M = \dot M^\perp$ and $\dot M_{22} = 0$.
\end{Lemma}
\begin{proof}
A variation of (I) together with Lemma \ref{Lem:usinh} gives
\[ \left( \sinh ((n-1)s) (M' M^{-1})^{\cdot} \right)' = 0. \]
Hence 
\[ (M' M^{-1})^{\cdot} = \frac{1}{\sinh (n-1)s} P \]
for some constant matrix $P$.
Moreover, observe that
\begin{equation} \label{eq:Mdotprime} \left( M' M^{-1} \right)^{\cdot} = \dot M' M^{-1} - M' M^{-1} \dot M M^{-1} = M ( M^{-1} \dot M )' M^{-1}.  
\end{equation}
So $(M' M^{-1})^\cdot$ is mapped to $-(M' M^{-1})^\cdot$ under the transformation $\perp$ combined with $s \to -s$.
Since $\sinh (n-1)s $ is odd, this implies that $P = P^\perp$, i.e. that $P$ is of block form
\[ P = \left( \begin{matrix}
P_{22} & 0 & \cdots & 0 \\
0 & P_{33} & \cdots & P_{3n} \\
\vdots & \vdots & \ddots & \vdots \\
0 & P_{n3} & \cdots & P_{nn}
\end{matrix} \right). \] 

Since by (\ref{eq:Mdotprime}) the lower block of $(M' M^{-1})^{\cdot}$ stays bounded for $s \to 0$, we find that the lower block of $(\sinh (n-1)s)^{-1} P$ must also stay bounded.
Hence, the lower block of $P$ must be zero.
Furthermore by (\ref{eq:Mdotprime}) and $\tr M^{-1} \dot M = \tr \dot M M^{-1} \equiv 2 \dot A$, we find $\tr (M' M^{-1} )^{\cdot} = 0$ and hence $\tr P = 0$.
So $P = 0$ and we conclude using (\ref{eq:Mdotprime}) again that $\dot M M^{-1} = Q$ for some constant matrix $Q$.

Now again since the problem is symmetric with respect to the transformation $\perp$ combined with $s \to -s$ and constant functions are even, we conclude that $Q = Q^\perp$.
Since moreover $M_{22}(s) = Q_{22} V(r(s))$, we conclude by smoothness at the origin that $Q_{22} = 0$.
\end{proof}

We can now summarize the discussion above:
Returning to the old coordinates $(r, \theta, x_3, \ldots, x_n)$, we have proven so far that $h$ takes the following form
\[ h =  - \DIV^* \xi + r^2 \sum_{i,j=3}^n u_{ij} \dd x_i \dd x_j. \]
So $h_{22} = - (\DIV^* \xi)_{22}$.
By the equations from the proof of Lemma \ref{Lem:xi}, we conclude from the decay of $h$ that $V^{1/2} \xi_1 (r) < C r^{-n+1.1}$ hence $r^{-2} |(\DIV^* \xi)_{ii}| < C r^{-n+1.1}$ for $i \geq 3$.
Together with $|h| < C r^{-n+1.1}$ this implies $u_{ij} = 0$ and thus $k \equiv 0$.
Hence by Lemma \ref{Lem:xi} we have $h \equiv 0$.
This concludes the proof of Proposition \ref{Prop:EinstvarBH}.

\end{document}